\documentclass[a4paper,12pt]{amsart}

\usepackage[latin1]{inputenc}
\usepackage{url}
\usepackage[french]{babel}
\usepackage{amsmath,amssymb,amsfonts,amsthm,amsfonts}

\def\]{\textup{\mbox{]\hspace{-.15em}]}}}
\def\[{\textup{\mbox{[\hspace{-.15em}[}}}

\newtheorem{lemme}[subsection]{Lemme}
\newtheorem{prop}[subsection]{Proposition}

\newtheorem{cor}[subsection]{Corollaire}

\theoremstyle{definition}

\newtheorem{remarque}[subsection]{Remarque}

\newenvironment{pf}
        {\medskip\noindent {\it Preuve --- \ }}
        {\hfill\nobreak $\Box$ \par\bigbreak}

\input xy
\xyoption{all}

\title{Quelques courbes de Hecke se plongent dans l'espace de Colmez}

\newcommand{\Q}{\mathbb{Q}}
\newcommand{\Qb}{\overline{\Q}}

\newcommand{\Gp}{{\rm G}_p}
\newcommand{\Z}{\mathbb Z}
\newcommand{\F}{\mathbb F}
\newcommand{\W}{\mathcal W}
\newcommand{\Wp}{{\rm W}_{{\Q_{p}}}}

\newcommand{\C}{\mathcal C}
\newcommand{\CC}{\mathbb C}
\newcommand{\G}{\mathbb G}
\newcommand{\Gm}{{\mathbb G}_{m}}

\newcommand{\ps}{\par \smallskip}

\newcommand{\phig}{(\varphi,\Gamma)}
\newcommand{\Gal}{{\rm Gal}}

\newcommand{\Qpb}{\overline{\Q}_{p}}
\newcommand{\Fpb}{\overline{{\mathbb F}}_p}
\newcommand{\Qp}{\Q_{p}}
\newcommand{\Dr}{{\mathrm D}_{{\rm rig}}}
\newcommand{\Dc}{{\mathrm D}_{{\rm cris}}}
\newcommand{\Fil}{{\rm Fil}}

\newcommand{\Ro}{{\mathcal R}}
\newcommand{\N}{\mathbb N}
\newcommand{\Ext}{{\mathrm Ext}}
\newcommand{\Fb}{\overline{\mathbb{F}}}
\newcommand{\rec}{{\rm rec}}
\newcommand{\Rep}{{\rm Rep}}
\newcommand{\figa}{{\rm FG}}
\newcommand{\rk}{{\rm rg}}
\newcommand{\GL}{{\rm GL}}
\newcommand{\HH}{{\mathcal H}}
\newcommand{\ZZ}{\mathcal Z}
\newcommand{\OO}{\mathcal O}
\newcommand{\isomo}{\overset{\sim}{\rightarrow}}

\newcommand{\ZT}{\widetilde{\ZZ(U_{p})}}
\newcommand{\LF}{\mathcal L_{\rm F}}
\newcommand{\D}{\rm D}
\newcommand{\Spec}{{\rm Spec}}

\newcommand{\Tr}{{\rm Tr}}

\author[G.~Chenevier]{Ga\"etan Chenevier}
\email{chenevie@math.univ-paris13.fr}
\address{C.N.R.S. \\ Laboratoire Analyse G\'eom\'etrie et Applications \\
Institut Galil\'ee \\ Universit\'e Paris 13 \\
99 Av. J.-B. Cl\'ement \\
93430 Villetaneuse \\ France}
\thanks{L'auteur est financ\'e par le C.N.R.S.\, Courriel : \og \url{chenevie@math.univ-paris13.fr}\fg\,.}

\begin{document}

\maketitle

{\sc Abstract : } Let $p$ be a prime, $\C$ the $p$-adic Eigencurve (with tame level $1$) and $\ZT$ the blow-up of the Fredholm hypersurface of the $U_p$-operator at the {\it special} points. We show that for $p = 2, 3, 5$ and $7$, the natural map  $$\C \longrightarrow \ZT$$ is a rigid-analytic isomorphism.\ps
\bigskip
\bigskip
\section*{Introduction}
Soient $p$ un nombre premier et $\C$ la courbe de Hecke $p$-adique de $\GL_{2}$ de niveau mod\'er\'e $1$ \,\,--\, \,{\it the Eigencurve}\,\,--\,\,construite par 
Coleman et Mazur dans \cite{cm}\footnote{Il y a deux constructions de la courbe de Hecke {\it loc. cit.} , nous utiliserons dans cet article uniquement celle qu'ils notent $D$.}. Cette courbe rigide analytique $p$-adique $\C$ (r\'eduite, s\'epar\'ee) admet par construction un morphisme localement fini et plat 
$$\kappa: \C \longrightarrow \W,$$  o\`u $\W$ d\'esigne l'espace de caract\`eres $p$-adiques \og des poids\fg\, de $\Z_{p}^{*}$. 

Soit $\HH$ la $\Z$-alg\`ebre de polyn\^{o}mes abstraite engendr\'ee par les symboles $U_p$ et $T_l$ pour tout $l \neq p$ premier, ce qui fournit des $T_n \in \HH$ pour tout entier $n \geq 1$ par les formules usuelles\footnote{Formellement, $(1-U_p\,p^{-s})^{-1}\prod_{l \neq p} (1-T_l\,l^{-s}+l^{\kappa-1-2s})^{-1} = \sum_{n\geq 1} T_n\,n^{-s}$.}. \`A chaque $T\in \HH$ 
est associ\'ee une fonction analytique born\'ee par $1$ sur $\C$, que l'on d\'esignera du m\^{e}me nom. 
Pour chaque $x \in \C$, il existe une unique forme modulaire surconvergente $f_{x}$ propre pour $\HH$, de pente finie, et normalis\'ee\footnote{Strictement, il faut exclure le point Eisenstein ordinaire 
$p$-adique de poids $0$ pour que cela soit vrai, c'est sans importance ici.}, telle que $T(f_x)=T(x)f_x$ pour tout $T \in \HH$, {\it i.e.} telle que pour tout $n \geq 1$ on ait $$a_{n}(f_{x})=T_{n}(x),$$ 
o\`u $a_{n}(g)$ d\'esigne le $n^{ieme}$ coefficient de $g$ pour son $q$-d\'eveloppement usuel \`a l'infini. 
Enfin, chaque telle forme appara\^{i}t ainsi ; cela inclut notamment toutes les formes modulaires usuelles de niveau $\Gamma_1(p)$, et l'ensemble des points $x$ sous-jacents \`a celles-l\`a est en fait Zariski-dense dans $\C$. \par
Il se trouve que par construction $U_{p}$ est une fonction inversible sur $\C$, ce qui permet de consid\'erer le morphisme
$$U_{p}^{-1}: \C \longrightarrow \W \times \G_{m},$$
qui s'av\`ere \^{e}tre fini (mais non surjectif \'evidemment) par construction. Toujours par construction, 
son image est la courbe r\'eduite $$\ZZ(U_{p}) \subset \W \times \G_{m}$$
d\'efinie comme \'etant la courbe spectrale de l'op\'erateur $U_{p}$. Rappelons que suivant Coleman,
l'op\'erateur $U_{p}$ agit de mani\`ere compacte 
sur la famille orthonormale d'espaces de Banach sur $\W$ constitu\'ee des espaces de formes modulaires surconvergentes, et que (toujours suivant Coleman) on peut donc consid\'erer sa s\'erie caract\'eristique $$\det(1-TU_{p}) \in \OO(\W \times \G_{m}),$$
dont $\ZZ(U_{{p}})$ est par d\'efinition le lieu des z\'eros. \ps \medskip

{\bf Definition:} {\it Un point $(\kappa,\lambda^{-1}) \in \W \times \G_m$ sera dit {\it sp\'ecial}, si $\kappa=\chi^k=(x \mapsto x^k)$ pour un certain entier $k\geq 2$ et si de plus $\lambda^{2}=p^{k-2}$. Plus g\'en\'eralement, on dira que $x \in \C$ est sp\'ecial si $(\kappa(x),U_p^{-1}(x))$ l'est. }\par \smallskip

Le th\'eor\`eme de classicit\'e de Coleman~\cite{colemanclassicite} assure que si $x \in \C$ est sp\'ecial, $f_x$ est une forme modulaire classique de poids $k$ et niveau $\Gamma_0(p)$. La relation $a_p(f_x)^2=p^{k-2}$ entra\^ine, pour des raisons de poids, que cette forme est n\'ecessairement $p$-nouvelle. La r\'eciproque \'etant \'evidente, il vient que {\it $x$ est sp\'ecial si, et seulement si, $f_x$ est une forme modulaire parabolique classique de niveau $\Gamma_0(p)$ et $p$-nouvelle}. Pour faire court, nous appellerons tout simplement {\it sp\'eciale} une telle forme\footnote{Notons que cela \'equivaut encore \`a demander que $f_x$ engendre une repr\'esentation lisse {\it sp\'eciale} de $\GL_2(\Q_p)$ au sens usuel, d'o\`u la terminologie adopt\'ee ici.}.  \par \medskip

Inspir\'es par l'article de Colmez \cite{colmeztri}, 
consid\'erons l'\'eclatement $\widetilde{\W \times \G_m}$ de $\W \times \G_m$ aux points sp\'eciaux, ainsi que
$$\ZT \subset \widetilde{\W \times \G_m}$$ 
le transform\'e strict de $\ZZ(U_{p})$ par cet \'eclatement. Si $x \in \C$ est sp\'ecial, il se trouve que $\kappa: \C \rightarrow \W$ est \'etale en $x$ (cf. Lemme~\ref{specialetale}), et en particulier que {\it $x$ est un point lisse de $\C$}. On en d\'eduit que l'application naturelle $\C \longrightarrow \ZZ(U_{p})$ 
se rel\`eve en un morphisme rigide analytique 
\begin{equation}\label{mainmap}\C \longrightarrow \ZT.
\end{equation}

On se propose dans cet article de d\'emontrer le r\'esultat suivant : \ps \medskip

{\bf Th\'eor\`eme}\,:\, {\it Si $p = 2, 3, 5$ ou $7$, le morphisme (\ref{mainmap}) est un isomorphisme. \par \smallskip

De plus, en chaque point sp\'ecial $(\chi^k,\lambda^{-1})$, $\ZZ(U_{p})$ est un diviseur \`a croisements normaux. Ses branches sont en bijection naturelle avec les formes modulaires sp\'eciales propres pour $\HH$ de poids $k$ et de valeur propre $\lambda$ sous $U_p$, la direction d'une branche \'etant naturellement reli\'ee au $L$-invariant de Fontaine de la forme associ\'ee.}
\par \medskip \medskip

Ce th\'eor\`{e}me a les cons\'equences ponctuelles amusantes suivantes. \ps \ps

{\bf Corollaire }\,:\, {\it Supposons encore $p = 2, 3, 5$ ou $7$. Si $\kappa \in \W$, notons $M_\kappa^\dagger$ l'espace des formes modulaires $p$-adiques surconvergentes de poids $\kappa$ et niveau mod\'er\'e $1$.
\begin{itemize}
\item[(i)] Soit $f \in M_\kappa^\dagger$ une forme propre pour l'op\'erateur de Hecke $U_p$ de valeur propre $a_p \neq 0$. Si $(\kappa,a_p^{-1})$ n'est pas sp\'ecial, alors $f$ est propre pour tous les op\'erateurs de Hecke dans $\HH$. De plus, si $f=\sum_{n\geq 0} a_n(f) q^n$ est normalis\'ee alors pour tout $n\geq 1$  $$a_{n}(f) \in \Q_{p}(a_p,\kappa).$$
\item[(ii)] Soient $f$ et $g$ deux formes modulaires sp\'eciales propres pour $\HH$, de m\^{e}me poids et normalis\'ees. Alors $f=g$ si et seulement si $a_p(f)=a_p(g)$ et si elles ont m\^{e}me $L$-invariant $\LF$ de Fontaine\footnote{Ainsi que nous l'a fait remarqu\'e H. Hida, cet \'enonc\'e implique l'existence de beaucoup de formes modulaires sp\'eciales en $p$ dont l'invariant $\LF$ est non nul. En effet, le nombre de formes propres $p$-nouvelles de niveau $\Gamma_0(p)$ cro\^it asymptotiquement lin\'eairement avec le poids $k$ (pair), et pour chaque poids $k$ donn\'e il n'y a que deux possibilit\'es pour $a_p$, \`a savoir $\pm p^{k/2-1}$.} en $p$. De plus, pour tout $n\geq 1$, $$a_{n}(f) \in \Q_{p}(a_p(f),\,\LF(f)\,).$$
\item[(iii)] Si $\lambda \neq 0$ est une valeur propre de $U_p$ sur $M_\kappa^\dagger$, alors l'espace propre de $U_p$ associ\'e \`a $\lambda$ est de dimension $1$, \`a moins que $(\kappa,\lambda^{-1})$ ne soit sp\'ecial, auquel cas il est de dimension \'egale au nombre de formes sp\'eciales de tels poids $\kappa$ et $U_p$-valeur propre $\lambda$.\end{itemize}}\ps \bigskip

Le premier ingr\'edient dans la d\'emonstration du th\'eor\`eme consiste en une \'etude locale-globale des repr\'esentations galoisiennes $p$-adiques attach\'ees aux formes modulaires surconvergentes de niveau mod\'er\'e $1$, que nous effectuons au \S\ref{etapelocglob}. Pr\'ecis\'ement, nous montrons que pour $p\leq 7$, ces repr\'esentations galoisiennes $p$-adiques globales sont d\'etermin\'ees par un sous-groupe de d\'ecomposition en $p$, dans ce sens qu'elles ont m\^eme images que leur restriction \`a un tel sous-groupe (Prop.~\ref{resp}). Lorsque $p=2$, cette propri\'et\'e nous avait \'et\'e communiqu\'ee par M. Emerton, ainsi que l'id\'ee que $\C$ devrait du coup pouvoir s'\'etudier \og purement localement\fg\, dans ce cas.
En fait, ces valeurs de $p$ sont les seules pour lesquelles cette propri\'et\'e est satisfaite. Les arguments ici sont \'el\'ementaires et ont peu de liens avec les formes modulaires surconvergentes : ils reposent essentiellement sur la th\'eorie du corps de classes et sur la nullit\'e du genre de $X_1(p)$ et du groupe des classes de $\Q(\mu_p)$ pour ces $p$.\par
  \par\smallskip

Le second ingr\'edient consiste en l'\'etude ponctuelle et infinit\'esimale de la famille $p$-adique de repr\'esentations de $\Gal(\Qpb/\Qp)$ port\'ee par $\C$ ; cela occupe les sections $2$ \`a $5$, le nombre premier $p$ y est quelconque ainsi que le niveau mod\'er\'e $N$ auxiliaire. Pour chaque point $x \in \C$ et chaque \'epaississement artinien $\Spec(A) \hookrightarrow \C$ de $x$ dans $\C$, nous disposons d'une repr\'esentation galoisienne $$\rho_A: \Gal(\Qpb/\Qp) \longrightarrow \GL_2(A)$$ (canonique \`a quelques exceptions pr\`es) que nous voulons d\'ecrire aussi pr\'ecis\'ement que possible. Nous disposons pour cela de r\'esultats de Kisin \cite{kisin} assurant le prolongement analytique de certaines p\'eriodes cristallines h\'erit\'ees par Zariski-densit\'e de celles des points classiques. Ainsi que l'avait remarqu\'e Colmez \cite{colmeztri} ces r\'esultats de Kisin s'interpr\^etent efficac\'ement, du moins ponctuellement, en terme de propri\'et\'es de triangulation des $\phig$-modules associ\'es sur l'anneau de Robba. Ce point de vue avait alors \'et\'e repris et \'etendu au cadre infinit\'esimal (et \`a la dimension sup\'erieure) dans \cite[Ch. 2]{BCh}, et c'est celui-l\`a aussi que nous adoptons ici. 

Aussi, les sections \ref{fgmodules}, \ref{rescolmez} et \ref{resbch} sont consacr\'ees aux $\phig$-modules g\'en\'eraux. Nous rappelons au \S\ref{rescolmez} les r\'esultats de Colmez concernant les $\phig$-modules de rang $1$ et leurs extensions, ainsi que son crit\`ere de triangulation. Dans la section~\S\ref{resbch}, reprenant des arguments de \cite{BCh}, nous en donnons les analogues infinit\'esimaux. Dans la section \S\ref{applicationC}, nous explicitons l'application \`a la courbe de Hecke $\C$. Le point crucial, r\'esultant des propri\'et\'es de triangulation et des calculs d'extensions pr\'ec\'edents, est qu'\`a certaines exceptions \'evidentes pr\`es nous d\'ecrivons compl\`etement la repr\'esentation $\rho_A$ plus haut (Prop. \ref{description}, Prop. \ref{descriptioninf}) et ce {\it uniquement en terme du couple $(U_p,\kappa)$ sur $\Spec(A)$}. Mentionnons tout de m\^eme ici une subtilit\'e concernant les points (non classiques) $x$ de poids-caract\`ere $\kappa=x^{k}\chi(x)$ ($\chi$ d'ordre fini) pour $k\geq 2$ et tels que $v(U_p(x)) > k-1$ : leurs \'epaississements infinit\'esimaux ont des propri\'et\'es triangulines plus restreintes. Pr\'ecis\'ement, pour un tel $x$, $\rho_A$ est $A$-trianguline si, et seulement si, $\Spec(A)$ est dans la fibre de $\kappa$ en $x$. \ps\ps

La section~\ref{preuves} contient la preuve du th\'eor\`eme principal. Nous commen\c{c}ons par v\'erifier que $$\C \longrightarrow \ZT$$ est injectif sur les points ferm\'es. Un point $x$ de $\C$ \'etant uniquement d\'etermin\'e par le couple $(\rho_x,U_p(x))$ o\`u $\rho_x$ est la repr\'esentation galoisienne globale attach\'ee \`a $f_x$, il faut voir que ce dernier couple est uniquement d\'etermin\'e par la donn\'ee \og locale\fg\, de $(\kappa(x),U_p(x))$. En dehors de certains points exceptionnels trait\'es \`a part, cela r\'esulte alors des \S\ref{etapelocglob} et \S\ref{applicationC}. 
Nous v\'erifions ensuite que $\C \longrightarrow \ZT$ est une immersion ferm\'ee, {\it i.e.} s\'epare les vecteurs tangents. L'argument est alors similaire \`a l'argument ponctuel sauf que nous utilisons cette fois les propri\'et\'es triangulines infinit\'esimales (en fait sur les $A=k(x)[\epsilon]/(\epsilon^2)$) des $\rho_A$ donn\'ees au \S\ref{applicationC}. 
Nous \'etudions enfin les divers points exceptionnels manquants. Pour les premiers $p$ consid\'er\'es, le lieu ordinaire cuspidal de $\C$ est vide\footnote{C'est certainement attendu puisque la repr\'esentation globale est d\'etermin\'ee par un groupe de 
d\'ecomposition en $p$.}, et il est possible de d\'eterminer explicitement les points non ordinaires de $\C$ dont la repr\'esentation locale associ\'ee est r\'eductible : ce sont les points Eisenstein critiques, ils sont \'etales sur $\W$ car les $p$ ci-dessus sont r\'eguliers. L'analyse des points sp\'eciaux utilise un autre th\'eor\`eme de 
Colmez \cite{colmezLinv} reliant l'invariant $L$ local avec la d\'eriv\'ee de $U_p$ par rapport au poids. 

 \ps \medskip

Terminons par mentionner que, du moins lorsque $p=2$, la g\'eom\'etrie de la courbe $\C$ a \'et\'e \'etudi\'ee par de nombreux auteurs (Emerton, Smithline, Buzzard, Kilford, Jacobs, Calegari, Loeffler, voir \cite{BCa} et \cite{BK}), et ce par calculs explicites de l'action de la correspondance de Hecke $U_2$ sur certaines r\'egions de $X_1(2)(\CC_2)$. 
Les r\'esultats que nous obtenons sont compl\'ementaires. Il est par exemple int\'eressant de les comparer avec ceux de \cite{BK} dans lequel les auteurs d\'ecrivent l'image r\'eciproque dans $\C$ du bord $$\{\kappa \in \W, 1/8<|\kappa(5)-1|<1\}$$ de $\W$ : cette couronne ext\'erieure est en fait la plus grande \'evitant tous les points sp\'eciaux, en lesquels les valeurs propres sp\'eciales de $U_2$ ont en g\'en\'eral de la mutliplicit\'e. \ps

\bigskip
{\sc Remerciements.} L'auteur remercie Jo\"el Bella\"iche, Kevin Buzzard, Pierre Colmez, Matthew Emerton et Barry Mazur pour des discussions utiles, ainsi que les organisateurs et acteurs du semestre {\it Eigenvarieties} de Harvard university pendant lequel l'essentiel de cet article a \'et\'e \'ecrit. Nous remercions de plus le referee pour ses remarques pertinentes.

\section{Repr\'esentations $p$-adiques globales de dimension $2$ d\'etermin\'ees par un groupe de d\'ecomposition en $p$}\label{etapelocglob}

\subsection{D\'efinitions} Soient $G$ et $H$ deux groupes munis d'un morphisme $\iota : H \longrightarrow G$. Pour toute repr\'esentation de $G$, et plus g\'en\'eralement tout morphisme de groupes $\rho : G \longrightarrow \Gamma$, on dira que $\rho$ est {\it d\'etermin\'ee par sa restriction \`a $H$} si $\rho$ et $\rho \circ \iota$ ont m\^eme image. Il est \'evidemment \'equivalent de demander que l'application naturelle $H \longrightarrow G/{\rm Ker}\, \rho$ est surjective.\ps 

Si $K$ est un corps de nombres, on notera $G_K$ le groupe de Galois absolu de $K$ et $G_{K,S}$ le groupe de Galois d'une extension alg\'ebrique maximale de $K$ non ramifi\'ee hors de $S$. Pour toute place $v$ de $K$ on dispose d'une classe de conjugaison canonique de morphismes $$\Gal(\overline{K}_v/K_v) \longrightarrow G_{K,S}.$$ On dira alors qu'une repr\'esentation de $G_{K,S}$ est d\'etermin\'ee par un groupe de d\'e\-com\-po\-si\-tion (resp. d'inertie) en $p$, si elle est d\'etermin\'ee par sa restriction \`a $\Gal(\overline{K}_v/K_v)$ (resp. au groupe d'inertie de ce dernier), ce qui ne d\'epend pas du morphisme choisi. 

\subsection{Structure de certains groupes de Galois}\label{calccyclo} 
Fixons $p$ un nombre premier et $K$ un corps de nombres dans lequel $p$ ne se 
d\'ecompose pas. Soient $P$ l'unique premier de $K$ divisant $p$, et $K(p)$ la pro-$p$ extension galoisienne maximale de $K$ qui est 
non ramifi\'ee hors de $P$ et des places archim\'ediennes. Le lemme suivant est bien connu.

\begin{lemme}\label{regular} $P$ est totalement ramifi\'e dans $K(p)$ si, et seulement si, le groupe 
des classes strictes d'id\'eaux de $\OO_K$ est d'ordre premier \`a $p$.
\end{lemme}

\begin{pf} D'apr\`es le calcul classique du sous-groupe de Frattini d'un pro-$p$-groupe, un sous-groupe ferm\'e $I$ d'un pro-$p$-groupe $Q$ qui se surjecte dans le quotient
$$Q/\overline{[Q,Q]\langle Q^p \rangle}$$
est en fait \'egal \`a $Q$. Appliquons ceci \`a  $Q:=\Gal(K(p)/K)$ et \`a un sous-groupe d'inertie $I$ au dessus de $P$. Le sous-corps de $K(p)$ fix\'e par l'image de $I$ dans le quotient de Frattini \'etant la $(p,p,\cdots,p)$-extension ab\'elienne maximale de $K$ non ramifi\'ee \`a toutes les places finies, on conclut par la 
th\'eorie du corps de classes.
\end{pf}

\begin{cor}\label{corpN} Si $p$ est r\'egulier, le morphisme $$G_\Q \longrightarrow \Gal(\Q(\mu_p)(p)/\Q)$$ est d\'etermin\'e par un groupe d'inertie en $p$. En particulier, toute repr\'esentation continue de $G_{\Q,\{\infty,p\}}$, dont la restriction \`a $G_{\Q(\mu_{p})}$ est d'image pro-$p$, est d\'etermin\'ee par un groupe d'inertie en $p$.
\end{cor}

\subsection{Repr\'esentations galoisiennes en petits conducteur et caract\'eristique} 

\begin{lemme}\label{genre1} Le genre de $X_1(M)$ est nul si, et seulement si, $1\leq M \leq 10$ ou $M=12$.
\end{lemme}

\begin{pf} C'est bien connu, nous laissons le lecteur le d\'eduire par exemple de \cite[Thm. 4.2.5, 4.2.9, 4.2.11]{Mi}. Pour simplifier les calculs on peut remarquer que si $X_1(M)$ est de genre $0$, alors $X_1(d)$ l'est aussi pour tout $d$ divisant $M$. Si $M>3$ est premier, alors le genre de $X_1(M)$ est $1+(M-1)(M-11)/24$.
\end{pf}

Soit $p$ un nombre premier. Noter que d'apr\`es le lemme ci-dessus, $X_1(p)$ est de genre $0$ si et seulement si $p=2, 3, 5$ ou $7$. 

\begin{lemme} \label{serre} Supposons que $p = 2, 3, 5$ ou $7$. Alors toute repr\'esentation continue, semisimple, impaire $$G_{\Q,\{\infty,p\}} \longrightarrow \GL_{2}(\Fb_{p})$$ est r\'eductible. En particulier, elle se factorise par $\Gal(\Q(\mu_{p})/\Q)$.
\end{lemme}

\begin{pf} Supposons tout d'abord que cette repr\'esentation est attach\'ee \`a une forme modulaire propre 
$f$ de poids $k\geq 2$ et niveau $\Gamma_1(p)$. D'apr\`es un argument bien connu, on peut supposer $k=2$. Mais pour les $p$ de l'\'enonc\'e $X_1(p)$ est de genre $0$, de sorte que $f$ est une s\'erie d'Eisenstein, ce que l'on voulait. Ce simple fait sera bien entendu suffisant pour notre application \`a la courbe de Hecke. Le cas g\'en\'eral r\'esulte de la 
d\'emonstration r\'ecente de la conjecture de Serre en niveau $1$ par Khare\footnote{Quand $p=2$ (resp. $p=3$), c'est un th\'eor\`eme de Tate (resp. Serre).}.
\end{pf}
\newcommand{\rhob}{{\overline{\rho}}}
Soit $A$ une $\Q_p$-alg\`ebre finie, locale, de corps r\'esiduel $L:=A/m_A$, munie de sa topologie naturelle de $\Q_p$-alg\`ebre de Banach. Pour $p=2$, le r\'esultat suivant \'etait connu de M. Emerton.

\begin{prop} \label{resp} Soient $p = 2, 3, 5$ ou $7$, et $$\rho: G_{\Q,\{\infty,p\}} \longrightarrow \GL_2(A)$$ une repr\'esentation continue, impaire si $p>2$, alors : \begin{itemize}
\item[(i)] la restriction de $\rho$ \`a $\Gal(\Qb/\Q(\mu_{p}))$ est d'image pro-$p$,
\item[(ii)] $\rho$ est d\'etermin\'ee par sa restriction \`a un groupe d'inertie en $p$,
\item[(iii)] si $\rho': G_{\Q,\{\infty,p\}} \longrightarrow \GL_2(A)$ est une autre repr\'esentation ayant les m\^emes propri\'et\'es que $\rho$, alors $\rho$ et $\rho'$ sont isomorphes, si et seulement si, elles le sont apr\`es restriction \`a un groupe d'inertie en $p$.
\end{itemize}
\end{prop}

\begin{pf} Comme nos premiers $p$ sont r\'eguliers, (i) entra\^ine (ii) et (iii) d'apr\`es le Corollaire~\ref{corpN}. D\'emontrons (i). \par

Consid\'erons plus  g\'en\'eralement une repr\'esentation continue $$\rho : G \longrightarrow \GL_n(A),$$ 
$G$ \'etant un groupe topologique compact quelconque. Si $\rho \bmod m_A$ est d'image pro-$p$, il en va de m\^eme pour $\rho$ : cela vient de ce que les quotients successifs des groupes multiplicatifs $1+m^i_A M_n(A)$, pour $i\geq 1$, sont isomorphes au groupe additif $L^{n^2}$, dont tous les sous-groupes compacts sont pro-$p$. Supposons alors que $A=L$, puis par un choix de r\'eseau que $\rho(G) \subset \GL_n(\OO_L)$ ; par un raisonnement semblable au pr\'ec\'edent, $\rho$ est d'image pro-$p$ si son image dans $\GL_n(\Fpb)$ est un $p$-groupe. Revenant \`a notre probl\`eme, on conclut par le Lemme \ref{serre}.
\end{pf}

\begin{prop}\label{respnec} La liste ci-dessus est exhaustive : si $p \geq 11$, il existe une re\-pr\'e\-sen\-ta\-tion galoisienne impaire 
$$\rho: G_{\Q,\{\infty, p\}} \longrightarrow \GL_{2}(\Qb_{p}),$$ qui n'est pas d\'etermin\'ee par un groupe de d\'e\-com\-po\-si\-tion en $p$. 
\end{prop}

\begin{pf} En effet, consid\'erons la repr\'esentation $\rhob$ r\'esiduelle (semisimple) mod $p$ attach\'ee \`a la 
forme modulaire $\Delta$ de Ramanujan. Un r\'esultat de Swinnerton-Dyer \cite[Cor. Thm. 4]{SD} montre que 
$p=2, 3, 5, 7, 23$, et $691$ sont les seuls premiers exceptionnels de $\Delta$, {\it i.e.} tels que l'image de 
$\rhob$ ne contienne pas ${\rm SL}_{2}(\F_{p})$. Comme ce dernier n'est pas r\'esoluble, 
$\rho$ n'est pas d\'etermin\'ee par son groupe de d\'ecomposition en $p$ pour ces $p$. 
On \'elimine alors de m\^eme $p=23$ et $691$ en consid\'erant l'unique forme parabolique propre $\Delta_{16}$ de poids 
$16$ et niveau $1$, pour laquelle Swinnerton-Dyer montre que $23$ et $691$ ne sont pas exceptionnels.
\end{pf}

\section{$\phig$-modules sur l'anneau de Robba}\label{fgmodules}

\subsection{Notations et conventions}\label{notconvfgmod} Soient $p$ un nombre premier, $\Q_p$ le corps des nombres $p$-adiques, 
$\Qb_p$ une cl\^oture alg\'ebrique de $\Q_p$ munie de sa 
valuation telle que $v(p)=1$, $\CC_p$ la compl\'etion 
de $\Qb_p$, et $$\Gp=\Gal(\Qb_p/\Q_p)$$ 
le groupe de Galois absolu de $\Q_p$. 
On note aussi $\Wp \subset \Gp$ le groupe de Weil de $\Q_p$. 
Normalisons l'isomorphisme 
$$\rec_{p}: \Q_p^* \longrightarrow \Wp^{\rm ab}$$ de la th\'eorie du corps de classes locale de sorte que Frobenius g\'eom\'etriques et uniformisantes se correspondent. \ps

Soit $L$ une extension finie de $\Q_{p}$.  L'{\it anneau de Robba} \`a coefficients dans $L$
est la $L$-alg\`ebre $\Ro_L$ des germes de fonctions analytiques $$f(z)=\sum_{n\in \Z} a_n (z-1)^n, \, \, 
a_n \in L$$
convergeantes sur une couronne de $\CC_p$ de la forme $r_f \leq |z-1| < 1$, $r_f$ d\'ependant de $f$. D'apr\`es des r\'esultats de Lazard, cette $L$-alg\`ebre est un domaine de B\'ezout. 
Elle admet de plus une topologie localement convexe naturelle (limite inductive 
de limites projectives d'alg\`ebres de Banach sur $L$). 
On d\'efinit une action de $\varphi$ et $\Gamma:=\Z_p^*$ sur $\Ro_L$ par les formules
$$\varphi(f)(z)=f(z^p), \qquad \, \, \, \gamma(f)(z)=f(z^{\gamma}), \, \gamma \in \Gamma.$$
\noindent L'action de $\Gamma$ ainsi obtenue est continue et commute avec celle de $\varphi$. 
Un \'el\'ement important de $\Ro_L$ est 
$$t=\log(z):=\sum_{n\geq 1} (-1)^{n+1}\frac{(z-1)^n}{n}.$$
Cet \'el\'ement converge sur toute la boule unit\'e ouverte 
$\{z \in \CC_p, |z-1|<1\}$, et s'annule \`a l'ordre $1$ aux racines de l'unit\'e de cette boule. 
Il v\'erifie de plus $$\varphi(t)=pt, \, \, \, \, \, \gamma(t)=\gamma t,\, \forall 
\gamma \in \Z_p^{*}.$$ 
\par
Par {\it $L$-repr\'esentation} $V$ de $\Gp$ nous entendrons toujours une repr\'esentation continue, 
$L$-lin\'eaire, sur un $L$-espace vectoriel de dimension finie. Nous noterons $\Rep$ la cat\'egorie tensorielle de ces repr\'esentations. Nous conviendrons que le caract\`ere cyclotomique de 
$\Gp$, {\it i.e} $\Q_{p}(1)={\rm T}_p(\mu_{p^{\infty}})\otimes \Q_p$, a pour poids de Hodge-Tate $-1$ et polyn\^{o}me de Sen $X+1$. 

\subsection{$\phig$-modules} 
Un $\phig$-module (sous-entendu sur $\Ro_L$) est un $\Ro_L$-module libre de type fini $D$ muni d'actions $\Ro_L$-semi-lin\'eaires 
de $\varphi$ et $\Gamma$ telles que $\Ro_L\varphi(D)=D$, $\Gamma$ agit contin\^ument pour la topologie 
$\Ro_L$-lin\'eaire de $D$, et $\varphi$ et $\Gamma$ commutent. Les $\phig$-modules 
forment une $\otimes$-cat\'egorie de mani\`ere naturelle. 
Des travaux de Kedlaya \cite{kedlaya} d\'efinissent une th\'eorie de pentes pour les $\varphi$-modules $M$ libres de rang fini sur $\Ro_L$, 
$\varphi$-semi-lin\'eaires, et satisfaisant $\Ro_L\varphi(M)=M$. Un 
$\phig$-module est alors dit {\it \'etale} si son $\varphi$-module sous-jacent n'admet que des pentes nulles.
Des travaux de Fontaine \cite{fontaine} et Cherbonnier-Colmez 
permettent de construire un foncteur de la cat\'egorie $\Rep$ 
des $L$-repr\'esentations de dimension finie et continues de $\Gp$ vers la cat\'egorie $\figa$ des $\phig$-modules, que nous noterons $V \mapsto \Dr(V)$. 

\begin{prop}\label{equcat} (Fontaine, Cherbonnier-Colmez, Kedlaya, \cite[prop. 2.7]{colmeztri}) Le foncteur $\Dr$ induit une $\otimes$-\'equivalence de cat\'egories 
entre $\Rep$ et la sous-cat\'egorie pleine $\figa_{et}$ 
des objets \'etales de $\figa$. De plus, $\dim_L(V)=\rk_{\Ro_L}(\Dr(V))$ pour toute $L$-re\-pr\'e\-sen\-ta\-tion $V$.
\end{prop}

\begin{remarque}\label{remfg}\begin{itemize}
\item[(i)] Il faut prendre garde \`a ce que contrairement \`a sa sous-cat\'egorie $\figa_{et}$, $\figa$ 
ainsi d\'efinie n'est pas ab\'elienne, car par exemple $tD$ est toujours un sous-objet strict de $D$ 
n'admettant pas de quotient. Cependant, si l'on inverse $t$ ou si l'on rajoute les objets de torsion de 
pr\'esentation finie, elle devient ab\'elienne. Un des int\'er\^ets \`a \'elargir $\figa$ en $\figa_{et}$ 
vient de ce que certains objets \'etales se d\'evissent naturellement dans $\figa$ et non dans $\figa_{et}$.
\item[(ii)] Ainsi qu'il est observ\'e dans \cite[Lemme 2.2.5]{BCh}, les travaux de Kedlaya montrent que 
la sous-cat\'egorie $\figa_{et}$ est \'epaisse, {\it i.e} stable par extensions dans $\figa$.
\end{itemize}
\end{remarque}

\section{$\phig$-modules de rang $1$ et leurs extensions, d'apr\`es Colmez \cite{colmeztri}.}\label{rescolmez}

\subsection{Les $\phig$-modules $\Ro_L(\delta)$.} Soit $\delta: \Qp^* \longrightarrow L^*$ un caract\`ere continu, on lui associe suivant Colmez un $\phig$-module $\Ro_L(\delta)$ de rang $1$ sur $\Ro_L$ comme suit: son $\Ro_L$-module sous-jacent est $\Ro_L$ et on pose
$$\varphi(1)=\delta(p), \qquad \, \, \gamma(1)=\delta(\gamma), \gamma \in \Gamma.$$
La pente de Kedlaya de ce $\phig$-module est $v(\delta(p))$. Si cette pente est nulle, $\delta$ s'\'etend via l'isomorphisme du corps de classes local en un caract\`ere continu 
$\delta': \Gp \longrightarrow L^*$, auquel cas on verifie imm\'ediatement que $\Ro_L(\delta)=\Dr(\delta')$.\par

\begin{prop} \label{colrang1}\cite[prop. 4.2]{colmeztri} Tout $\phig$-module de rang $1$ est de la forme $\Ro_L(\delta)$ pour un unique $\delta$. 
\end{prop}

Cela d\'ecoule par exemple imm\'ediatement de l'\'equivalence de cat\'egories plus haut et du cas \'etale, mais Colmez en donne aussi une preuve directe. Colmez calcule aussi la $\phig$-cohomologie en degr\'es $0$ et $1$ des $\Ro_L(\delta)$, 
en terme du complexe \`a trois termes usuel. Ses r\'esultats ont \'et\'e \'etendus au degr\'e $2$, ainsi qu'en fait \`a tous les 
$\phig$-modules par Liu dans \cite[Thm. 0.2]{liu}. Soient $x: \Q_p^* \longrightarrow L^*$ l'inclusion naturelle, et $\chi: \Q_p^* \longrightarrow L^*$ le caract\`ere tel que $\chi(p)=1$ et $\chi_{|\Z_p^*}= x_{|\Z_p^*}$. 

\begin{prop} \label{colext}\cite[Thm. 3.9]{colmeztri} $\Ext_{\figa}(\Ro_L(\delta_2),\Ro_L(\delta_1))$ est de $L$-dimension $1$, 
\`a moins que $\delta_1\delta_2^{-1}$ ne soit de la forme $x^{-i}, \chi x^{i}$ pour $i\geq 0$ un entier, auquel cas il est de dimension $2$.
\end{prop}

Dans le cas $v(\delta_1\delta_2^{-1}(p))=0$, 
cela se d\'eduit du th\'eor\`eme de Tate, de l'\'equivalence de cat\'egorie plus haut, et du fait que $\figa_{et}$ est \'epaisse dans $\figa$. Cependant, le calcul de Colmez est ind\'ependant des r\'esultats de Tate. \par

\begin{remarque} 
Soit $D$ une extension non triviale de $\Ro_L(\delta_2)$ par $\Ro_L(\delta_1)$. 
Si $D$ est \'etale, le th\'eor\`eme de Kedlaya montre que $v(\delta_1\delta_2(p))=0$ et que $v(\delta_1(p))\geq 0$. 
R\'eciproquement, sous ces conditions, $D$ est \'etale d'apr\`es Kedlaya \`a moins qu'il ne 
contienne un sous-$\phig$-module $D'$ de rang $1$ (que l'on peut supposer satur\'e) de pente strictement n\'egative. 
Un tel module est n\'ecessairement en somme directe avec $\Ro_L(\delta_1)$ et $D$ est donc scind\'e sur $\Ro_L[1/t]$, 
ceci se produit de mani\`ere exceptionnelle (voir \cite{colmeztri} pour plus de d\'etails).
\end{remarque}
\subsection{Crit\`ere de triangulation de Colmez en rang $2$} Fixons $V$ une $L$-repr\'esentation continue de $\Gp$ de dimension $2$. Suivant Colmez, on dit que $V$ est {\it trianguline} si $\Dr(V)$ 
est une extension de $\Ro_L(\delta_2)$ par $\Ro_L(\delta_1)$ pour certains caract\`eres $\delta_1, \, \delta_2$. La proposition 
suivante est une version l\'eg\`erement raffin\'ee d'un r\'esultat de Colmez. 
\renewcommand{\Dc}{{\rm D}_{{\rm cris}}}
\begin{prop}\label{colmez} \cite[prop. 5.3]{colmeztri}\, \, $V$ est trianguline si, et seulement si, il existe un caractere continu $\psi: \Gp \longrightarrow L^*$ et un \'el\'ement 
$\lambda \in L^*$ tels que $\Dc(V\otimes \psi)^{\varphi=\lambda} \neq 0$. Sous cette derni\`ere condition, fixons un 
$v\neq 0 \in \Dc(V\otimes \psi)^{\varphi=\lambda} \neq 0$. Alors $\Dr(V)$ est une extension de $\Ro_L(\delta_2)$ par $\Ro_L(\delta_1)$, o\`u 
\begin{itemize}
\item[(i)] $\delta_1(p)=\lambda \psi^{-1}(p) p^{-s}$ et $\delta_{|\Gamma}=\psi^{-1}\chi^{-s}$, $s$ \'etant le plus petit entier tel que 
$v \notin \Fil^{s+1}(\Dc(V\otimes \psi))$.
\item[(ii)] $\delta_2=\delta_1^{-1}\det(V)$.
\end{itemize}
\end{prop}
\begin{pf} Le r\'esultat \'etant important pour la suite de cet article, nous redonnons la d\'emonstration de Colmez. 
Pr\'ecis\'ement, montrons la condition suffisante de triangulation, la r\'eciproque \'etant imm\'ediate. On peut supposer que $\psi=1$. 
On pose $D:=\Dr(V)$.\par
Un th\'eor\`eme de Berger \cite{berger} montre que le $L[\varphi]$-module $\Dc(V)$ s'identifie naturellement avec $(D[1/t])^{\Gamma}$. 
On dispose ainsi d'un \'el\'ement $v\neq 0 \in D[1/t]^{\Gamma}$. Quitte \`a multiplier $v$ par $t^{-s}$ pour un certain entier $s \in \Z$, 
$w:=t^{-s}v \in D$. L'id\'eal $I$ de $\Ro_L$ engendr\'e par les coefficients de $w$ dans une $\Ro_L$-base fix\'ee de $D$ est 
alors stable par $\varphi$ et $\Gamma$. Comme $\Ro_L$ est de B\'ezout et que $I$ est de type fini, $I$ est principal, et une analyse de son diviseur 
montre alors que $I=(t^i)$ pour un certain entier $i\geq 0$. Ainsi, quitte \`a diminuer $s$, on peut supposer que $I=\Ro_L$, i.e. que $\Ro_L w$ est 
facteur direct comme $\Ro_L$-module, par la propri\'et\'e de B\'ezout. \par
Consid\'erons le sous-$\phig$-module $D':=\Ro_L w$. Clairement, $D'\isomo \Ro_L(\delta_1)$ o\`u  $\delta_1$ satisfait les deux premi\`eres 
conditions de l'\'enonc\'e (i). La derni\`ere condition vient de la recette pour la filtration dans l'isomorphisme de Berger, et pr\'ecis\'ement du Lemme \cite[Lemme 2.4.2]{BCh}. \par
Mais $D/D'$ est alors libre de rang $1$, et donc de la forme $\Ro_L(\delta_2)$ par la classification. On conclut (ii) en prenant le d\'eterminant.
\end{pf}

\section{D\'eformations triangulines, d'apr\`es \cite{BCh}.}\label{resbch}
\newcommand{\deltab}{\overline{\delta}}
\subsection{Les $\phig$-modules $\Ro_A(\delta)$.} Dans cette section, $A$ est une $\Q_p$-alg\`ebre finie locale de corps r\'esiduel $A/m_A \isomo L$. On d\'efinit l'anneau de Robba \`a coefficients dans $A$ comme \'etant $$\Ro_A:=\Ro_{\Q_p} \otimes_{\Q_p} A.$$
Nous appellerons ici $\phig$-module sur $\Ro_A$ la donn\'ee d'un $\phig$-module sur $\Ro_{\Qp}$ muni d'une action de la $\Q_p$-alg\`ebre $A$ telle que le $\Ro_A$-module sous-jacent soit libre.\par

Soit $\delta: \Qp^* \longrightarrow A^*$ un caract\`ere continu. On d\'efinit $\Ro_A(\delta)$ de mani\`ere \'evidente comme plus haut (cf. \cite[\S 2.3.1]{BCh}). 
La remarque \ref{remfg} (ii) montre que $\Ro_A(\delta)$ est isocline de pente $v(\deltab(p))$, o\`u $$\deltab: \Qp^* \longrightarrow L^*$$ 
est la r\'eduction modulo $m_A$ de $\delta$. Comme plus haut, $\Ro_A(\delta)$ s'identifie au $\Dr$ du caract\`ere galoisien naturel dans 
le cas \'etale. En fait, toujours par la m\^{e}me remarque \ref{remfg} (ii), on montre comme plus haut la proposition suivante.

\begin{prop} \label{rg1def}\cite[prop. 2.3.1]{BCh} Tout $\phig$-module de rang $1$ sur $\Ro_A$ est de la forme $\Ro_A(\delta)$ pour un unique $\delta$.
\end{prop}

\renewcommand{\deltab}{\overline{\delta}}
\begin{prop} \label{extdef}Soient $\delta_1, \delta_2: \Qp^* \longrightarrow A^*$ deux caract\`eres continus. Supposons que $\deltab_1\deltab_2^{-1} \neq x^{-i}, \chi x^i$. Alors $\Ext_{\figa}(\Ro_A(\delta_2),\Ro_A(\delta_1))$ est un $A$-module libre de rang $1$.
\end{prop}

\begin{pf} On peut supposer $\delta_2=1$, et on note $\delta:=\delta_1$. Si $\deltab \neq x^{-i}, \chi x^i$ avec $i\geq 0$ entier,
alors d'apr\`es \cite[Thm. 3.9]{colmeztri} et \cite[Thm 0.2]{liu} on a $H^0(\Ro_L(\deltab))=H^2(\Ro_L(\chi\deltab^{-1})=0$, et $H^1(\Ro_L(\deltab))$ est de $L$-dimension $1$. Un d\'evissage 
imm\'ediat montre alors que $H^1(\Ro_A(\delta))$ est de $\Q_p$-dimension $\dim_{\Qp} A$. \par
Prouvons que ce $A$-module est libre de rang $1$ par r\'ecurrence sur la longueur de $A$. 
C'est imm\'ediat si $A$ est un corps, on peut donc supposer qu'il existe un id\'eal $I \subset m_A$ de longueur $1$. 
On regarde alors la suite exacte de $\phig$-modules:
$$0 \longrightarrow I\Ro_A(\delta) \longrightarrow \Ro_A(\delta) \longrightarrow \Ro_A(\delta \otimes A/I) \longrightarrow 0, $$
o\`u l'on a $I\Ro_{A}(\delta)\isomo\Ro_L(\deltab)$. Cette suite reste exacte apr\`es passage aux $H^1$ par le paragraphe pr\'ec\'edent. Ainsi, on voit que $I H^1(\Ro_A(\delta)) = H^1(I\Ro_A(\delta))$, et donc que comme $A$-module on a 
$$H^1(\Ro_A(\delta))\otimes_A A/I \simeq H^1(\Ro_{A/I}(\delta \otimes A/I)).$$
On conclut par induction et Nakayama que $H^1(\Ro_A(\delta))$ est $A$-monog\`ene, puis libre de rang $1$.
\end{pf}

\begin{remarque} \label{remdef}Supposons que $\delta: \Qp^* \longrightarrow A^*$, $A=L[\epsilon]/(\epsilon^2)$, est non constant et satisfait $\deltab=x^{-i}$ ou 
$\deltab=\chi x^i$. Dans ce cas, on montre facilement que $H^1(\Ro_A(\delta))$ est isomorphe \`a $L \oplus A$ comme $A$-module.
\end{remarque}

\subsection{Crit\`ere de triangulation artinien en rang $2$.} Fixons maintenant une
re\-pr\'e\-sen\-ta\-tion $A$-lin\'eaire continue $V_A$ de $\Gp$ qui est libre de rang $2$ comme $A$-module, on pose 
$D_A:=\Dr(V_A)$. Alors $D_A$ libre de rang $2$ sur $\Ro_A$ (\cite[Lemme 2.2.7]{BCh}), c'est donc un $\phig$-module sur $\Ro_A$ au sens plus haut. 
On dira que $V_A$ est {\it $A$-trianguline} si $D_A$ est une extension de $\Ro_A(\delta_1)$ par $\Ro_A(\delta_2)$ 
pour certains caract\`eres continus $\delta_1, \delta_2: \Q_p^* \longrightarrow A^*$. On pose $\overline{V}:=V_A/m_A V_A$. 
Le crit\`ere de triangulation est un peu plus subtile dans le cas artinien que dans le cas ponctuel (Prop. \ref{colmez}):

\begin{prop} \label{tridef} \cite[Lemme 2.5.2]{BCh} Supposons que pour un $\lambda \in A^*$, il existe un $A$-module libre de rang $1$ $A.v \subset \Dc(V_A)^{\varphi=\lambda}$, et 
que la filtration de Hodge induite sur $A.v$ n'ait qu'un seul saut, disons $s\in \Z$. Alors $D_A$ est une extension de $\Ro_A(\delta_2)$ par $\Ro_A(\delta_1)$ 
o\`u $\delta_1(p)=\lambda p^{-s}$ et  $(\delta_1)_{|\Gamma}=\chi^{-s}$, et $\delta_2=\det(V_A)\delta_1^{-1}$. En particulier, $V_A$ est $A$-trianguline.
\end{prop}

\begin{pf} D'apr\`es la recette de Berger pour la filtration, pr\'ecis\'ement par \cite[Lemme 2.4.2]{BCh}, 
on sait que $v \in t^sD_A$. De plus, 
le poids de $\overline{v} \subset Av$ \'etant aussi \'egal \`a $s$
par hypoth\`ese, $\Ro_{L}\overline{v}t^{{-s}}$ est satur\'e dans $\Dr(\overline{V})$. D'apr\`es \cite[Lemme 2.2.3]{BCh}, il vient que $\Ro_{A}t^{{-s}}v$ est satur\'e dans $D_{A}$ et on conclut comme dans la preuve de la Proposition
\ref{colmez}.
\end{pf}

\begin{remarque}\label{remdeftri} \begin{itemize}
\item[(i)] Pour qu'un $\phig$-module de la forme $\Ro_{A}(\delta)$ soit cristallin, il faut et suffit que $\delta_{{|\Z_p^{*}}}$ soit constant 
et \'egal \`a $x^{i}$ pour un certain entier $i\in \Z$.
\item[(ii)] Supposons que $A.v \subset \Dc(V_A)^{\varphi=\lambda}$ soit libre de rang $1$ et que le poids de 
$L\overline{v} \subset \Dc(V)^{\varphi=\overline{\lambda}}$ soit le plus petit poids de Hodge-Tate (entier) de $V$, disons $k$. 
Dans ce cas, la filtration de Hodge sur $A.v$ a pour unique saut $k$. En effet, une r\'ecurrence sur la longueur de $A$ montre que $\Fil^{k+1}(Av)=0$.
\end{itemize}
\end{remarque}

\section{Application \`a la courbe de Hecke de $GL_{2}$.}\label{applicationC}

Dans cette section nous donnons une application des trois sections pr\'ec\'edentes \`a la famille de repr\'esentations galoisiennes port\'ee par les courbes de Hecke.
\newcommand{\GG}{{\rm Gal}(\Qb/\Q)_{\{\infty,pN\}}}
Fixons $p$ un nombre premier, $N \geq 1$ un entier premier \`a $p$, et $\C$ la courbe de Hecke $p$-adique de niveau mod\'er\'e 
$\Gamma_{1}(N)$ comme dans l'introduction (\`a ceci pr\`es que nous autorisons ici un niveau mod\'er\'e $N$ quelconque). Cette courbe d\'epend d'un choix d'alg\`ebre de Hecke commutative $\HH$, et il sera suffisant pour nos besoins de supposer que $\HH$ contient $U_p$ et les $T_l$ pour $(l,Np)=1$. Les points $x \in \C$ sont encore en bijection naturelle avec les syst\`emes de valeurs propres de $\HH$ sur les espaces de formes modulaires $p$-adiques surconvergentes de niveau mod\'er\'e $N$ et de pente finie, et on d\'esignera par $f_x \neq 0$ n'importe laquelle des formes propres de poids $\kappa(x)$ ayant le syst\`eme de valeurs propres attach\'e \`a $x$. \ps

D'apr\`es l'existence connue des repr\'esentations galoisiennes
attach\'ees aux formes modulaires classiques (Eichler-Shimura, Igusa, Deligne), on peut d\'efinir un pseudocaract\`ere continu, de dimension $2$,
 $$\Tr : \GG \rightarrow \OO(\C)$$
dont l'\'evaluation en tout point $z \in \C$ param\'etrant une forme $f_{z}$ qui est classique est la trace de la repr\'esentation 
galoisienne $p$-adique attach\'ee \`a cette forme.
En retour, on dispose {\it pour tout} $x \in \C$ d'une unique repr\'esentation galoisienne continue, semi-simple,
$$\rho_{x}: \GG \rightarrow \GL_{2}(k(x)),$$
telle que la trace d'un Frobenius g\'eom\'etrique en un premier $l$ ne divisant pas $Np$ est $T_{l}(x)$
(il n'y a pas d'obstruction de Schur \`a cause de la conjugaison complexe, de sorte que $\rho_{x}$ est d\'efinie sur le corps r\'esiduel $k(x)$).
Les propri\'et\'es connues en $p$ des repr\'esentations attach\'ees aux formes modulaires et des travaux de Sen, Kisin, 
montrent que le polyn\^{o}me de Sen de $\rho_{x}$ (sous-entendu restreinte \`a un groupe de d\'ecomposition en $p$) 
est $$X(X-d\kappa(x)+1).$$ Ici, $d\kappa(x)$ signifie la d\'eriv\'ee en $1$ de $\kappa(x)$ vu comme caract\`ere continu 
de $\Z_{p}^{*}$. De plus, on a (cf.  \cite[Thm 6.3]{kisin})
\begin{equation} \label{mk} \Dc^{+}(\rho_{x})^{\varphi=U_{p}(x)} \neq 0.\end{equation}
Enfin d\'esignons par $\W_{N}$ l'espace rigide analytique sur $\Q_{p}$ param\'etrant les caract\`eres continus $p$-adiques 
de $$\Z_{p}^{*}\times (\Z/N\Z)^{*}=\GG^{\rm ab}.$$
Ces caract\`eres peuvent \^{e}tre vue de mani\`ere naturelle comme des caract\`eres continus de $\Q_{p}^{*}$ en \'etendant l'identit\'e de $\Z_{p}^{*}$ en un morphisme
$$\Q_{p}^{*} \rightarrow \Z_{p}^{*}\times (\Z/N\Z)^{*},\, \, \, \, p \mapsto (1,p).$$
Par la th\'eorie du corps de classes locale, on peut donc aussi voir les caract\`eres de $\Z_{p}^{*}\times (\Z/N\Z)^{*}$ comme des caract\`eres de $\Gp$.

\subsection{Propri\'et\'es triangulines ponctuelles}

Fixons $L/\Q_{p}$ une extension finie et $x \in \C(L)$, correspondant donc \`a une forme modulaire 
$f_{x} \neq 0$ surconvergente de niveau mod\'er\'e $\Gamma_{1}(N)$ propre pour $\HH$ et de pente finie, 
poids-caract\`ere $$\kappa_{N}(x):=(\kappa(x),\varepsilon(x)) \in \W_{N}(L).$$ Par construction, 
$\det(\rho_{x})=\kappa_{N}^{-1}\chi$.
Posons $k:=d\kappa(x)$, $\lambda:=U_{p}(x)$, $V=\rho_{x}$ 
restreinte \`a $\Gp$, $D=\Dr(V)$ le $\phig$-module sur $\Ro_{L}$ attach\'e \`a $V$.

\begin{prop} \label{description} Il existe deux caract\`eres continus $\delta_{1}, \delta_{2}: \Q_{p}^{*} \longrightarrow L^{*}$ et 
une suite exacte de $\phig$-modules sur $\Ro_{L}$ 
\begin{equation}\label{eqtrig} 0 \longrightarrow \Ro_{L}( \delta_{1} ) \longrightarrow D \longrightarrow \Ro_{L}(\delta_{2}) \longrightarrow 0
\end{equation}
o\`u $\delta_{2}=\det(V)\delta_{1}^{-1}$ et $\delta_{1}$ est comme suit. D\'efinissons un caract\`ere $\delta: \Q_{p}^{*} \rightarrow L^{*}$ par $\delta(p)=\lambda$ et $\delta$ trivial sur $\Z_{p}^{*}$.
\begin{itemize}
\item[(i)] Si $k \notin \Z$, ou $k\in -\N$, ou $k$ est un entier $\geq 1$ et $v(\lambda) < k-1$, alors $\delta_{1}=\delta$. De plus, dans ce cas:
\begin{itemize}
\item Si $v(\lambda) = 0$, alors $V$ contient une $L$-droite non ramifi\'ee sur laquelle le Frobenius g\'eom\'etrique agit par multiplication par $\lambda$. 
\item Si $v(\lambda)\neq 0$, alors l'extension~(\ref{eqtrig}) n'est pas scind\'ee et est unique, \`a moins que $k\geq 2$ soit un entier, 
$\kappa(x)=\chi^{k}$, $\lambda^{2}\varepsilon(p)=p^{k-2}$, auquel cas $V$ est semi-stable et uniquement d\'etermin\'ee par la donn\'ee suppl\'ementaire de son invariant $L$ de Fontaine.
\end{itemize}
\item[(ii)] Si $k$ est un entier $\geq 1$ et $v(\lambda)>k-1$, alors $\delta_{1}=x^{1-k}\delta$ et l'extension (\ref{eqtrig}) est non scind\'e et unique. De plus $D$ n'est pas De Rham.
\item[(iii)] Si $k$ est un entier $\geq 1$ et $v(\lambda)=k-1 \geq 0$, alors $V$ est r\'eductible. Soit $\psi$ le caract\`ere non ramifi\'e de 
$\Gal(\Qpb/\Qp)$ sur lequel le Frobenius g\'eom\'etrique agit par multiplication par $p^{k-1}/\lambda\varepsilon_{N}(p)$. Alors,
\begin{itemize}
\item soit $\Fil^{k-1}(\Dc(V))^{\varphi=\lambda}\neq 0$, auquel cas $V$ contient 
$L(1-k) \otimes_{L} \psi^{-1}$.
\item soit $\Fil^{k-1}(\Dc(V))^{\varphi=\lambda}=0$, auquel cas $k>1$ et $V$ est une rep\'esentation ordinaire non scind\'ee 
contenant $\psi$ comme sous-repr\'esentation non ramifi\'ee.
\end{itemize}
\end{itemize}
\end{prop}

\begin{pf} On obtient une triangulation explicite de $D$ en utilisant (\ref{mk}) et la Proposition \ref{colmez}, modulo la d\'etermination de l'entier $s \in \Z$. 
D\'eterminons ce $s$. Par (\ref{mk}), on sait que $s\geq 0$. 
De plus, $s$ \'etant un saut de la filtration de Hodge, c'est un poids de Hodge-Tate de $V$. Ainsi, $s=0$ sauf peut-\^{e}tre si $k$ est un entier $>1$, auquel cas $s$ peut \^{e}tre \'egal \`a $k-1$.
Supposons maintenant que $k>1$ est un entier. Notons que d'apr\`es \cite[prop. 2.3.4]{BCh}, on sait que si $s=0$ alors $V$ est de De Rham (donc potentiellement semistable).
\par
Supposons tout d'abord que $v(\lambda)<k-1$. Par un th\'eor\`eme de Coleman, on sait que $f_{x}$ est classique, et donc le polygone de Hodge de $\Dc(V)^{\varphi=\lambda}$ est en dessous de son polygone de Newton, ce qui montre que $s=0$.
Une autre mani\`ere de raisonner est de dire que si $s=k-1$, $D$ contient un sous-objet de pente $v(\lambda)+1-k<0$, ce qui contredit le th\'eor\`eme de Kedlaya. \par
Si maintenant $v(\lambda)>k-1$, alors $V$ n'est pas De Rham, i.e. non potentiellement semi-stable, 
car (\ref{mk}) contredirait la faible admissibilit\'e de ${\rm D}_{\rm pst}(V)$. Ainsi, $s=k-1$.\par
Enfin, si $v(\lambda)=k-1$, les deux possibilit\'es peuvent se produire (et se produisent en pratique). Si $s=0$, $V$ est de De Rham, et l'inspection des polygones de Hodge et Newton du foncteur de Fontaine montre que 
$V$ est ordinaire non scind\'ee. Si $s=k-1$, il vient que $V$ est r\'eductible est contient une droite sur laquelle $G_p$ agit par un twist non ramifi\'e de $L(1-k)$, de quotient un caract\`ere non ramifi\'e.\par
Les assertions d'unicit\'e  d\'ecoulent alors de la Proposition \ref{colext}.
\end{pf}

\begin{remarque} Les points $x$ de $\C$ tels que $\rho_{x}=\chi_{1}+\chi_{2}$ est (globallement) r\'eductible 
peuvent \^{e}tre d\'etermin\'es tout comme dans \cite[\S 4.2]{lissite}, ils sont tous attach\'es \`a des s\'eries d'Eisenstein classiques si $v(U_{p}(x))\neq 0$ (resp. dans les familles Eisenstein ordinaires sinon). 
\end{remarque}
\subsection{Propri\'et\'es triangulines infinit\'esimales}\label{propinf}

Fixons $x \in \C$ comme dans la section pr\'e\-c\'edente, et fixons de plus un \'epaississement infinit\'esimal 
d'anneau $A$ de $x$ dans $\C$. L'anneau $A$ est donc local artinien de corps r\'esiduel $L$. Supposons que : \begin{itemize}
\item[-] $\rho_x$ est irr\'eductible (hypoth\`ese globale), 
\item[-] si $v(\lambda)=0$, alors\footnote{En fait, nos arguments fonctionnent verbatim sous l'hypoth\`ese plus faible \og ${\rho_x}|_{\Gp}$ est ramifi\'ee\fg.} $k \neq 1$.  \end{itemize}
Dans ce cas, par la premi\`ere hypoth\`ese, 
il existe une repr\'esentation canonique 
$$\rho_{A}: \GG \rightarrow \GL_{2}(A)$$
de trace la compos\'ee $$\GG \overset{\Tr}{\longrightarrow} \OO(\C) \overset{{\rm can}}{\longrightarrow} A.$$ 
Notons encore $V_{A}$ la restriction de $\rho_{A}$ \`a un groupe de d\'ecomposition en $p$, $V=V_{A}/m_{A}V_A$. 
\begin{prop} \label{descriptioninf} Soit $\delta: \Q_{p}^{*} \rightarrow A^{*}$ le caract\`ere trivial sur $\Z_{p}^{*}$ et tel que $\delta(p)=U_{p}$.
\begin{itemize}
\item[(i)] Supposons que l'on est dans le cas (i) de la Proposition \ref{description}. Alors $V_{A}$ est une d\'eformation trianguline de $V$. Pr\'ecis\'ement, $\Dr(V_{A})$ est une extension de $\Ro_{A}(\delta_{2})$ par $\Ro_{A}(\delta_{1})$ relevant la triangulation~(\ref{eqtrig}) {\rm loc. cit.} de $\Dr(V)$, o\`u $\delta_{1}=\delta$ et $\delta_{1} \delta_{2}=\kappa_{N}^{-1}\chi$. Si de plus $x$ n'est pas sp\'ecial et $v(\lambda)\neq 0$, alors une telle extension est unique.
\item[(ii)] Supposons que l'on est dans le cas (ii) de la Proposition \ref{description}. Alors 
$V_{A}$ est $A$-trianguline si, et seulement si, ${\rm Spec}(A)$ est dans la fibre de $\kappa$ 
au dessus du point ferm\'e $\kappa(x) \in \W$. Si tel est le cas, $\Dr(V_{A})$ est une extension de $\Ro_{A}(\delta_{2})$ par 
$\Ro_{A}(\delta_{1})$ relevant la triangulation de $\Dr(V)$, o\`u $\delta_{1}=x^{1-k}\delta$ et $\delta_{1}
\delta_{2}=\kappa_{N}^{-1}\chi$. Une telle extension est unique.
\end{itemize}
\end{prop}

\begin{pf} La th\'eorie de Sen montre que le polynome de Sen de $V_{A}$ est $X(X-d\kappa+1) \in A[X]$. Les travaux de Kisin montrent que 
$$\Dc^{+}(V_{A})^{\varphi=U_{p}}$$
est libre de rang $1$ sur $A$ d\`es que $\Dc^{+}(V)^{\varphi=\lambda}$ 
est de $L$-dimension $1$ (en ces termes, voir par exemple \cite[Thm. 3.3.2]{BCh}). 
C'est toujours le cas \`a moins que $V$ ne soit cristalline de poids $(0,0)$ ({\it i.e.} non ramifi\'ee), ce qui est exclu par hypoth\`ese sur $x$. \par
Supposons que l'on est dans le cas (i), la Proposition \ref{tridef} conclut \`a l'aide de la remarque \ref{remdeftri} . \par
Supposons que l'on est dans le cas (ii), et que ${\rm Spec}(A)$ est inclus dans la fibre de $\kappa$ en $x$. Dans ce cas, la th\'eorie de Sen montre que 
$$(V_{A}\otimes_{\Q_{p}}\mathbb C_{p}(k-1))^{\Gp}$$
\noindent est libre de rang $1$ sur $A$. Les th\'eor\`emes de Tate sur la cohomologie de $\mathbb C_{p}(i)$ assurent alors 
($k-1$ \'etant le plus grand poids de Hodge-Tate de $V$)
que $$\Fil^{k-1}\D_{\rm DR}(V_{A})$$
est libre de rang $1$ sur $A$. Mais comme d'autre part $\D_{\rm DR}(V)$ est de $L$-dimension $1$, cela montre que
$$\Dc(V_{A})^{\varphi=U_{p}}=\Fil^{k-1}(\D_{\rm DR}(V_{A})).$$
Ainsi, la Proposition \ref{tridef} conclut la triangulation de $V_{A}$ dans ce cas. V\'erifions maintenant la r\'eciproque, {\it i.e.} supposons que l'on dispose d'une suite exacte
$$0 \longrightarrow \Ro_{A}(\delta_{1})\longrightarrow \Dc(V_{A}) \longrightarrow \Ro_{A}(\delta_{2}) \longrightarrow 0.$$
En r\'eduisant modulo l'id\'eal maximal de $A$, on obtient une triangulation de $V$ qui est n\'ecessairement\footnote{En effet, sinon $\Dr(V)[1/t]$ serait scind\'e, donc cristallin d'apr\`es la Proposition \ref{description}, ce qui est absurde car
$V$ n'est pas De Rham par la m\^{e}me proposition.} la triangulation donn\'ee par la Proposition \ref{description}. En particulier, le Frobenius cristallin de $\Dc(\Ro_{L}(\overline{\delta_{2}}))$ agit par multiplication par 
$$p^{k-1}/\lambda\varepsilon_{N}(p) \neq \lambda$$ 
(car $v(\lambda)>k-1$). Ainsi, il vient ais\'ement que $\Dc(\Ro_{A}(\delta_{1}))=\Dc(V_{A})$ est libre de rang $1$ sur $A$, et donc 
que $d\kappa-k+1$ est constant dans $A$ par la remarque \ref{remdeftri} (i).\par
	Les assertions d'unicit\'e dans (i) et (ii) d\'ecoulent alors de la Proposition~\ref{extdef}.
\end{pf}

\section{Preuve du th\'eor\`eme}\label{preuves}

Dans cette section, on prouve le th\'eor\`eme \'enonc\'e dans l'introduction. On commence par expliquer le fait suivant, admis dans l'introduction, qui vaut pour tout $p$ (mais $N=1$) : 

\begin{lemme}\label{specialetale} Le morphisme $\kappa$ est \'etale aux points sp\'eciaux.
\end{lemme}
\begin{pf} Soient $x \in \C$ un point sp\'ecial, $L:=k(x)$ le corps r\'esiduel de $x$, et $k:=d\kappa(x) \geq 2$. Comme $v(U_{p}(x))=k/2-1<k-1$, un th\'eor\`eme de Coleman \cite{colemanclassicite} assure que l'espace caract\'eristique de $f_{x}$ sous l'action de l'alg\`ebre de Hecke $\HH$ co\"incide avec son espace caract\'eristique dans l'espace des formes modulaires classiques de poids $k$ et niveau $\Gamma_0(p)$. Mais cet espace est de dimension $1$ sur $L$ par le principe du $q$-d\'eveloppement, car tous les op\'erateurs de Hecke (y compris $U_p$) y sont semisimples, ce qui conclut par \cite[Prop. 1]{lissite}.\end{pf}

Nous utiliserons librement les notations du \S\ref{applicationC}. Un point $x \in \C$ est dit {\it ordinaire} si $v(U_p(x))=0$, {\it r\'eductible} si $\rho_x$ l'est. Nous renvoyons \`a \cite[\S 4]{lissite} pour l'\'etude des points r\'eductibles de $\C$ : ce sont les points Eisenstein.

\begin{lemme} \label{petitlemme} Supposons $p \leq 7$.
\begin{itemize} 
\item[(i)] Le lieu ordinaire de $\C$ co\"incide avec la droite Eisenstein ordinaire.
\item[(ii)] Si $x \in \C$ est un point non ordinaire, $\rho_{x}$ est r\'eductible restreinte \`a un groupe de d\'ecomposition en $p$ si, et seulement si, $f_{x}$ est une s\'erie d'Eisenstein critique. En un tel point, $\C$ est lisse et $\kappa$ est \'etale. De plus, les points Eisenstein critiques sont les seuls points $x \in \C$ en lesquels on a une \'egalit\'e $v(U_{p}(x))=d\kappa(x)-1>0$. 
\end{itemize}
\end{lemme}

\begin{pf} Le premier point vient de la th\'eorie de Hida et de ce qu'il n'y a pas de forme modulaire parabolique ordinaire de poids $2$ sur $X_{1}(p)$ pour ces premiers, et sur $X_1(4)$ : ces courbes sont de genre $0$ par le Lemme~\ref{genre1}.\par
Le second point vient de ce que les points non ordinaires $x$ tels que $\rho_{x}$ est r\'eductible sont les points Eisenstein critiques, qui sont lisses et \'etales sur $\W$ pour ces premiers (cf. \cite[\S 4.1]{lissite}, car $p$ est r\'egulier), et que si $p\leq 7$
alors pour tout $x$, $\rho_{x}$ est d\'etermin\'ee par sa restriction \`a un groupe de d\'ecomposition en $p$ d'apr\`es la Proposition \ref{resp}.  Ainsi, $\rho_{x}$ est r\'eductible si, et seulement si, il l'est restreint \`a un groupe de d\'ecomposition en $p$. La derni\`ere assertion de (ii) d\'ecoule alors de la Proposition 
\ref{description} (iii).\par
\end{pf}

Commen\c{c}ons tout d'abord par un lemme sur le corps de d\'efinition des $\rho_{x}$, $x \in \C$. Rappelons que si $x$ est sp\'ecial, $\LF(x)$ est l'invariant $L$ de Fontaine de la repr\'esentation semistable non cristalline $V_x={\rho_x}|_{\Gal(\Qpb/\Qp)}$.

\begin{lemme}\label{rationalite} Pour tout $x\in \C$, $\rho_{x}$ est d\'efini sur $\Q_{p}(U_{p}(x),\kappa(x))$ si $x$ n'est pas 
un point sp\'ecial, sur $\Q_{p}(U_{p}(x),\LF(x))$ sinon.
\end{lemme}

\begin{pf} On peut supposer que $x$ n'est pas r\'eductible. D'apr\`es la Proposition \ref{description}, $\Dr(V_{x})$ est d\'efini sur le corps de l'\'enonc\'e, ainsi donc que $V_{x}$, puis $\rho_{x}$ par la Proposition \ref{resp} (ii). 
\end{pf}

D\'esignons par $$\pi: \C \longrightarrow \ZT,$$
le morphisme d\'efini dans l'introduction, gr\^ace au Lemme~\ref{specialetale}.

\begin{lemme} \label{injec} 
$\pi$ est bijectif sur les points (ferm\'es). 
\end{lemme}

\begin{pf} La surjectivit\'e de $\pi$ \'etant \'evidente par construction, v\'erifions son injectivit\'e.\ps
Soit $L$ une extension finie de $\Q_p$. Par construction de $\C$, deux points $x, y$ de $\C(L)$ sont \'egaux si, et seulement si, $U_p(x)=U_p(y)$ et $\rho_x \simeq \rho_{y}$. Par la Proposition \ref{resp} (iii), ce dernier isomorphisme \'equivaut \`a $V_{x} \simeq V_{y}$.

Soient $x, y \in \C(L)$ tels que $\pi(x)=\pi(y)$. En particulier, $U_{p}(x)=U_{p}(y)$, $\kappa(x)=\kappa(y)$ et il faut voir que $V_{x} \simeq V_{y}$. Posons $\lambda=U_{p}(x)$, $w=\kappa(x)$.
Comme $\kappa(x)=\kappa(y)$, on sait d\'ej\`a que $\det(V_{x})=\det(V_{y})$. \par
On peut supposer que $v(\lambda)>0$ car sinon $x$ et $y$ sont 
sur le lieu ordinaire, {\it i.e.} sur la droite Eisenstein ordinaire par le Lemme \ref{petitlemme}, 
auquel cas il est clair que $x=y$. 
Pour la m\^eme raison, d'apr\`es le point (ii) du Lemme \ref{petitlemme} on peut supposer que ni $x$, ni $y$, n'est Eisenstein critique, auquel cas le m\^{e}me lemme montre
que $v(\lambda) \neq dw -1$ et que $V_{x}$ et $V_{y}$ sont irr\'eductibles. Mezalor la Proposition \ref{description} 
(dont la disjonction des cas ne fait intervenir que $v(\lambda)$ et $w$) 
assure qu'\`a moins que $x$ et $y$ ne soient tous les deux sp\'eciaux, $V_{x} \simeq V_{y}$. \par
Supposons donc maintenant que 
$x$ et $y$ sont sp\'eciaux de m\^{e}me poids-caract\`ere $w$, et tels que $\lambda=U_{p}(x)=U_{p}(y)=\pm p^{k/2-1}$. Rappelons qu'en $x$ et $y$, $\kappa$ est \'etale par le Lemme \ref{specialetale}, 
de sorte que $\C$ est localement isomorphe \`a une boule ouverte dans $\W$ de centre $w$ : localement sur cette boule, $U_p$ est une fonction analytique de $\kappa$.
Le th\'eor\`eme principal de Colmez dans \cite{colmezLinv} assure que pour 
$z=x, y$, l'invariant de Fontaine $\LF(z)$ de $V_{z}$ vaut
\begin{equation}\label{formulecolmez}\LF(z)=-2\lambda^{-1}\left(\frac{\partial U_{p}}{\partial \kappa}\right)_{\kappa=w}= 2\lambda\left(\frac{\partial U_{p}^{-1}}{\partial \kappa}\right)_{\kappa=w}.\end{equation}
Ainsi, la Proposition \ref{description} (i) conclut que $V_{x} \simeq V_{y}$.
\end{pf}

\medskip

Puisque nous nous pla\c{c}ons en niveau mod\'er\'e $1$, rappelons que les courbes $\C$, $\ZZ(U_p)$, et donc $\ZT$, sont r\'eduites. Un morphisme entre $\Q_p$-espaces rigides r\'eduits qui est bijectif sur les points ferm\'es n'\'etant pas n\'ecessairement un isomorphisme\footnote{Penser
par exemple \`a la normalisation d'un cusp.}, il nous faut encore 
montrer que $\pi$ est une immersion ferm\'ee, ce qui concluera le 
premier point du th\'eor\`eme. Le second point d\'ecoule d\'ej\`a des Lemmes \ref{specialetale}, \ref{injec} et de la formule (\ref{formulecolmez}). \par \smallskip

Soit $z \in \C$ un point ferm\'e. Par le Lemme~\ref{rationalite}, $\pi$ induit un isomorphisme entre les corps r\'esiduels de $z$ et de $\pi(z)$, il ne reste donc qu'\`a montrer que l'application naturelle 
$$d\pi_{z}: T_{z}(\C) \longrightarrow T_{\pi(z)}(\ZT)$$
est injective. C'est imm\'ediat lorsque $\kappa$ est \'etale en $z$, 
car c'est alors vrai apr\`es composition par 
$T_{\pi(z)}(\ZT) \rightarrow T_{\kappa(z)}(\W)$. On peut donc supposer que $z$ n'est pas un point Eisenstein, {\it i.e.} par le Lemme \ref{petitlemme} (ii) que 
$\rho_z$ et $V_{z}:={\rho_{z}}_{|\Gp}$ sont irr\'eductibles, et aussi que $z$ n'est pas sp\'ecial.\ps
Alors que le Lemme~\ref{injec} reposait sur les r\'esultats de Colmez (\cite{colmezLinv},\cite{colmeztri}), le lemme suivant utilise la th\'eorie d\'evelopp\'ee dans \cite{BCh}.

\begin{lemme} \label{immersion} 
\begin{itemize}
\item[(i)] \`A moins que $k:=d\kappa(z)$ ne soit un entier $\geq 1$ et que $v(U_{p}(z))>k-1$,  $d\pi_{z}$ est injective.
\item[(ii)] Si $k:=d\kappa(z)$ est un entier $\geq 1$ et si $v(U_{p}(z))>k-1$, 
$d\pi_{z}$ est injective sur le sous-espace de $T_{z}(\C)$ correspondant aux vecteurs tangents 
dans la fibre de $\kappa: \C \longrightarrow \W$ au dessus de $\kappa(z)$.
\end{itemize}
Ainsi, dans tous les cas, $d\pi_z$ est injective sur le sous-espace d\'efini par $d\kappa_z=0$.
\end{lemme}

\begin{pf} 
Soient $A=k(z)[\epsilon]/(\epsilon^2)$, et $x, y$ deux $A$-points \'epaississements de $z$ dans $\C$ ayant m\^{e}me image dans 
$T_{\pi(z)}(\W \times \Gm)$, 
suppos\'es de plus dans la fibre de $\kappa$ vers $\W$ dans le cas (ii). Nous voulons montrer qu'ils sont \'egaux dans $T_{z}(\C)$. Rappelons que comme $\rho_{z}$ 
est irr\'eductible, il y a un sens \`a consid\'erer $\rho_{x}$ et $\rho_{y}$ comme dans la Proposition \ref{descriptioninf}, ainsi que $V_{x}$ et $V_{y}$ 
leur restriction respective \`a un groupe de d\'ecomposition en $p$.
Par le m\^{e}me argument que dans le Lemme \ref{injec}, il faut montrer que $V_{x} \simeq_A V_{y}$. Mais cela d\'ecoule de la Proposition \ref{descriptioninf} (noter que $z$ n'est pas sp\'ecial et 
que de plus on n'est jamais dans le cas (iii)
du Lemme \ref{description} \`a cause du Lemme \ref{petitlemme}).
\end{pf}

Terminons maintenant la d\'emonstration du th\'eor\`eme. Soit $z \in \C$, on dispose d'un diagramme commutatif naturel :

$$\xymatrix{ T_z(\C) \ar@{->}[dr]_{d\kappa_z} \ar@{->}[rr]^{d\pi_z} & & T_{\pi(z)}(\ZT) \ar@{->}[dl]^{d{\rm pr}_{1\,\pi(z)}} \\  & T_{\kappa(z)}(\W) & }$$
de sorte que ${\rm Ker}\,(d\pi_z) \subset {\rm Ker}\,(d\kappa_z)$, d'o\`u l'injectivit\'e de $d\pi_z$ par la derni\`ere assertion du Lemme~\ref{immersion}. $\square$
\ps \medskip
D\'emontrons maintenant le corollaire du th\'eor\`eme \'enonc\'e dans l'introduction. 

\begin{pf} Soit $M_\kappa^\dagger[\lambda] \subset M_\kappa^\dagger$ l'espace caract\'eristique de $U_p$ de valeur propre $\lambda$. Supposons d'abord $(\kappa,\lambda^{-1})$ non sp\'ecial. D'apr\`es le th\'eor\`eme principal, $\HH \otimes \Q_p(\kappa)$ et $\Q_p(\kappa)[U_p]$ ont m\^eme image dans ${\rm End}_{\Q_p(\kappa)}(M_\kappa^\dagger[\lambda])$. En particulier, un vecteur propre pour $U_p$ y est propre pour tous les $T \in \HH$. Par le \og principe du $q$-d\'eveloppement\fg\,, un tel vecteur propre est uniquement d\'etermin\'e \`a un scalair pr\`es par son syst\`eme de valeurs propres, du moins s'il est $p$-adiquement parabolique, ce que l'on peut supposer ici par le Lemme~\ref{petitlemme} (i). Cela montre le premier point de (i) et de (iii). Le second point de (i), ainsi que (ii) et le reste de (iii), d\'ecoulent des Lemmes~\ref{rationalite} et~\ref{injec}.
\end{pf}

\begin{remarque} \begin{itemize}
\item[(i)] Le th\'eor\`eme admet sans doute des g\'en\'eralisations partielles \`a d'autres niveaux mod\'er\'es $N$ (plut\^ot que $N=1$) et premiers $p$ ; nous n'avons pas cherch\'e \`a les incorporer par souci de clart\'e. Ainsi, si $N=2$ et $p=3, 5$, $N=3, 5$ et $p=2$, $N=4$ et $p=3$, la courbe modulaire $X_1(Np)$ est de genre $0$ et il existe une variante du th\'eor\`eme principal. On peut imaginer aussi d'autres variantes sp\'ecifiant certaines composantes connexes de $\W_N$, exploitant le fait que pour une poign\'ee suppl\'ementaire de couples ({\it niveau,caract\`ere}) $=(M,\chi)$ on a $S_2(\Gamma_0(M),\chi)=0$ (dans tous les cas, $M\leq 32$ comme on pourrait le v\'erifier par les formules de \cite{CO} et les tables de Stein~\cite{ST}). 

\item[(ii)] Une autre direction possible, disons encore en niveau mod\'er\'e $1$ pour simplifier, serait de fixer une repr\'esentation modulaire $$r: G_{\Q,\{\infty,p\}} \rightarrow \GL_{2}(\Fb_{p}),$$ par exemple irr\'eductible, soumise \`a la condition suppl\'ementaire que $$H^{1}(\Z[1/p],{\rm ad}(r)) \rightarrow H^{1}(\Q_{p},{\rm ad}(r))$$ est injective. On pourrait alors montrer que l'ouvert ferm\'e $\C(r) \subset \C$ correspondant \`a $r$ s'immerge encore dans l'\'eclatement de $\W \times \Gm$ aux points sp\'eciaux, du moins hors du lieu ordinaire.  
\end{itemize}
\end{remarque}

\end{document}